\documentclass{article}
\usepackage{amssymb}
\usepackage{amsmath}
\usepackage{latexsym}
\usepackage{mathrsfs}%花体
\usepackage{graphicx}%插图
\usepackage[all]{xy}%交换图
\newtheorem{Th}{Theorem}
\newtheorem{Lemma}[Th]{Lemma}
\newtheorem{Coro}[Th]{Corollary}

\title{A bijection proof of the Capelli's identity}
\author{XIONG Rui}
\begin{document}

\maketitle

\newcommand{\fracp}[2][]{\dfrac{\partial #1}{\partial #2}}
\def\S{\mathcal{S}}

\begin{abstract}
  In this article, a short combinatorial proof of the Capelli's identity is given.
  It also leads to an easy proof of the Capelli–-Cauchy–-Binet identity, a more general form of Capelli's identity.
  With the technique introduced, the Turnbull's identity can be proven by sum with signs.
  At the end of the article, the general Cayley identity is discussed as a dual version of Capelli's identity.
\end{abstract}

%\tableofcontents

\section{Introduction}

Let $\{x_{ij}:1\leq i,j\leq n\}$ be a set of indeterminants. Denote
$$D_{ij}=\sum_{k=1}^n x_{ik}\fracp{x_{jk}}. $$
The classical \emph{Capelli's identity} claimed
$$\left|\begin{matrix}
D_{11}+(n-1) & D_{12} & \cdots & D_{1n}\\
D_{21} & D_{22}+(n-2) & \cdots & D_{2n}\\
\vdots & \vdots & \ddots & \vdots\\
D_{n1} & D_{n2} & \cdots & D_{nn}
\end{matrix}\right|
=\left|\begin{matrix}
x_{11} & \cdots &x_{1n}\\
\vdots & \ddots & \vdots\\
x_{n1} &\cdots & x_{nn}
\end{matrix}\right|
\cdot
\left|\begin{matrix}
\fracp{x_{11}} & \cdots &\fracp{x_{1n}}\\
\vdots  & \ddots & \vdots\\
\fracp{x_{n1}} & \cdots & \fracp{x_{nn}}
\end{matrix}\right|. $$
The noncommutative determinant is expanded from left to right, say the left hand side is
$$\sum_{\sigma\in \S_n}(-1)^{\sigma} A_{\sigma(1)1}\cdots A_{\sigma(n)n},
\qquad A_{ij}=D_{ij}+\delta_{ij}(n-i). $$

The Capelli's identity was introduced more than one century ago \cite{capelli1887zuruckfuhrung}.
It turned out to be a useful tool in representation theory, see for example \cite{zbMATH00051906} appendix F for a glance.
It has several different proofs, including a combinatorial one, like \cite{Foata_1994}.
Different level of generalizations of it also vary, for example \cite{Caracciolo_2009} where
much more general results about noncommutative determinants are established.
Last, we discussed about Cayley identity, which is also classical see \cite{vivanti1890alcune}.
For its development, see for example \cite{Caracciolo_2013}.
The main idea is to generalize determinant to noncommutative case, or more specifically, the ring of differential operators.

In this article, a new proof through an algorithm is given.
The main method is by constructing a bijection between monomials and a certain set of configurations.
In some sense, the proof is just to cleverly count the correction to be commutative.
We can use the algorithm to get the Capelli–-Cauchy–-Binet identity.
To prove its variant, the Turnbull's identity, we should slightly use the trick of sum with signs.
%Last, we will discuss Cayley's identity.

%I am not sure whether the similar method can be used to prove Turnbull's identity, or Howe–-Umeda–-Kostant–-Sahi identity.

\section{The algorithm}

\bigbreak
\def\Fix{\operatorname{Fix}}
We call $(\sigma,\varphi)$ a \emph{Capelli configuration}
\begin{itemize}
  \item a permutation $\sigma\in \S_n$;
  \item assign each $i$ with $\sigma(i)=i$ a number $\varphi(i)\in \{1,\cdots, i\}$.
\end{itemize}
This can be illustrated by a diagram, for example,
%$$\begin{cases}
%\sigma=(1536), \\
%\varphi(2)=1, \\
%\varphi(4)=2, \\
%\varphi(6)=6;
%\end{cases}: \begin{array}{c}\xymatrix{\\
%1\ar@/^2pc/[rrrr]& 2\ar[l]& 3\ar@/^2pc/[rrrr]& 4 \ar@/^1pc/[ll] &5\ar@/^1pc/[ll]& 6\ar@(ul,ld)[] &7\ar@/^3pc/[llllll]\\&}
%\end{array}.$$
$$\begin{cases}
\sigma=(12583), \\
\varphi(4)=3, \\
\varphi(6)=6, \\
\varphi(7)=4;
\end{cases}:
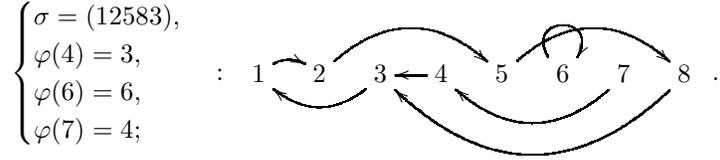
\begin{array}{c}\xymatrix@=1pc{\\
1\ar@/^0.5pc/[r]& 2\ar@/^1.5pc/[rrr]& 3\ar@/^1pc/[ll]& 4 \ar[l] &5\ar@/^1.5pc/[rrr]& 6\ar@(ul,ur)[]  &7\ar@/^1.5pc/[lll]&8\ar@/^2.5pc/[lllll]\\&}
\end{array}.$$

\bigbreak
Let $0\leq m\leq n$. Denote $\mathcal{C}^m$ all Capelli configurations $(\sigma,\varphi)$ such that
$$\left.\begin{array}{r}
\sigma(k)=k\\
\varphi(k)<k\end{array}\right\}\Longrightarrow k\geq m. $$
That is, $\varphi$ coincides with $\sigma$ over ${\Fix \sigma\cap \{1,\ldots,m-1\}}$. 
We also denote $\mathcal{C}=\mathcal{C}^{n+1}$, and note that
$$(\sigma,\varphi)\in \mathcal{C}\iff \varphi=\sigma|_{\Fix \sigma}. $$
Now, we define an algorithm
$$\Lambda: \quad \mathcal{C}^1
\stackrel{\Lambda^1}\longrightarrow\cdots \stackrel{\Lambda^{n-1}}\longrightarrow\mathcal{C}^n
\stackrel{\Lambda^{n}}\longrightarrow\mathcal{C}. $$
also $\Lambda=\Lambda^n\circ \cdots\circ \Lambda^1$ for convention.

\bigbreak
Let $(\sigma,\varphi)\in \mathcal{C}^m$ be a Capelli configuration.
When $\sigma(m)=m$ and $\varphi(m)<m$, then output $(\sigma',\varphi')$ such that
\begin{itemize}
  \item $\sigma'(\varphi(m))=m$, and $\sigma'(m)=\sigma(\varphi(m))$, and do not change other elements.
  \item $\varphi'$ is just restriction of $\varphi$ on $\Fix\sigma\setminus \{m\}$.
\end{itemize}
In other case, just output $(\sigma,\varphi)$ itself.

For the example above, if we run $\Lambda^4$, we will get
$$\begin{array}{c}
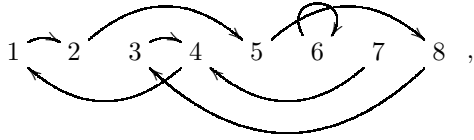
\xymatrix@=1pc{\\
1\ar@/^0.5pc/[r]& 2\ar@/^1.5pc/[rrr]& 3\ar@/^0.5pc/[r]& 4 \ar@/^1.5pc/[lll] &5\ar@/^1.5pc/[rrr]& 6\ar@(ul,ur)[]  &7\ar@/^1.5pc/[lll]&8\ar@/^2.5pc/[lllll]\\&}
\end{array},
$$
and then run $\Lambda^7$, we will get
$$\begin{array}{c}
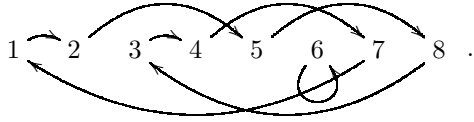
\xymatrix@=1pc{\\
1\ar@/^0.5pc/[r]& 2\ar@/^1.5pc/[rrr]& 3\ar@/^0.5pc/[r]& 4 \ar@/^1.5pc/[rrr] &5\ar@/^1.5pc/[rrr]& 6\ar@(dl,dr)[]  &7\ar@/^2pc/[llllll]&8\ar@/^2pc/[lllll]\\&}
\end{array}.
$$

\section{The construction}
Let $\{y_{ij}:1\leq i,j\leq n\}$ be another set of indeterminants commuting with $x_{ij}$'s.
Denote for a permutation $\sigma\in \S_n$
$$\Delta_{ij}=\sum_{k=1}^n x_{ik}y_{jk}.$$
If $G$ is a polynomial in $\{x_{ij}\}$, and $H$ in $\{y_{ij}\}$, and $F(x,y)=G(x)H(y)$, define
$$F(\partial)= G(x)H\bigg(\fracp{x}\bigg). $$
This extends to a linear map to all polynomials in $\{x_{ij}\}\cup \{y_{ij}\}$.
We will denote for $0\leq m\leq n$, and a Capelli configuration $(\sigma,\varphi)$
$$\Delta(\sigma,\varphi,m)
=\tilde{\Delta}_{\sigma(n)n}(\partial) \cdots \tilde{\Delta}_{\sigma(m)m}(\partial) \cdot
\big(\tilde{\Delta}_{\sigma(m-1),m-1}\cdots \tilde{\Delta}_{\sigma(1)1}\big)(\partial),$$
where
$$\tilde{\Delta}_{\sigma(i)i}= \begin{cases}
\Delta_{\sigma(i)i}, & \textrm{if $i\neq \sigma(i)$ or $\varphi(i)=i$},\\
1, & \textrm{otherwise}.
\end{cases}$$

\begin{Th}Under the notations above, for $(\sigma,\varphi)\in \mathcal{C}^{m+1}$, we have
$$(-1)^{\sigma}\Delta(\sigma,\varphi,m+1)=\sum_{\Lambda^m(\sigma' ,\varphi')=(\sigma,\varphi)}(-1)^{\sigma' }\Delta(\sigma' ,\varphi',m). $$
\end{Th}
Note that there are at most two of such $(\sigma' ,\varphi')$.
Firstly, $(\sigma,\varphi)$ itself such that $\Lambda^m(\sigma,\varphi)=(\sigma,\varphi)$.
By our algorithm,
that some other $(\sigma' ,\varphi')$ satisfies $\Lambda^m(\sigma' ,\varphi')=(\sigma,\varphi)$ only happens when $\sigma^{-1}(m)<m$ and
$\sigma' (\sigma^{-1}(m))=\sigma(m)$, and $\varphi'(m)=\sigma^{-1}(m)$.

When $\sigma^{-1}(m)\geq m$, then
$x_{m}$ commutes with $\fracp{x_{\sigma(i)}}$ for all $1\leq i\leq m-1$, so
$$\Delta(\sigma,\varphi,m+1)=\Delta(\sigma,\varphi,m). $$
When $i=\sigma^{-1}(m)< m$, then
$$\begin{array}{l}
\qquad
(\cdots)\bigg(x_{\sigma(m)}\fracp{x_{m}}\bigg)\cdot \bigg(x_{\sigma(1)}\cdots x_{\sigma(m-1)}\fracp{x_{1}}\cdots \fracp{x_{m-1}}\bigg)\\
=(\cdots)\bigg(x_{\sigma(m)}\cdot x_{\sigma(1)}\cdots x_{\sigma(m-1)}\cdot \fracp{x_{m}}\cdot \fracp{x_{1}}\cdots \fracp{x_{m}}\bigg)\\
\quad +
(\cdots)\bigg(x_{\sigma(m)}\cdot x_{\sigma(1)}\cdots \widehat{x_{\sigma(i)}}\cdots x_{\sigma(m-1)}\fracp{x_{1}}\cdots \fracp{x_{m-1}}\bigg), \\
=(\cdots)\bigg(x_{\sigma(1)}\cdots x_{\sigma(m-1)}\cdot \fracp{x_{m}}x_{\sigma(m)}\cdot \fracp{x_{1}}\cdots \fracp{x_{m}}\bigg)\\
\quad +
(\cdots)\bigg(x_{\sigma(1)}\cdots \underbrace{{x_{\sigma(i)}}}_{\textrm{exchange to }x_{\sigma(m)}}\cdots x_{\sigma(m-1)}\fracp{x_{1}}\cdots \fracp{x_{m-1}}\bigg). \\
\end{array}$$
%Note that the second term changing $x_{\sigma(i)}$ to $x_{\sigma(m)}$,
%which is exactly $\Delta(\sigma',\varphi',m)$.
That is,
$$\Delta(\sigma,\varphi,m)=\Delta(\sigma,\varphi,m+1)+\Delta(\sigma',\varphi',m), $$
where $(\sigma',\varphi')$ the only other Capelli configuration.
Clearly, $\sigma'$ differs from $\sigma$ by a swap.

\begin{Coro}\label{MainResultofCapelli}
Under the notations above, for $(\sigma,\varphi)\in \mathcal{C}$, we have
$$(-1)^{\sigma}\Delta(\sigma,\varphi,n+1)=\sum_{\Lambda(\sigma' ,\varphi')=(\sigma,\varphi)}(-1)^{\sigma' }\Delta(\sigma' ,\varphi',1). $$
\end{Coro}

\bigbreak
Since when $(\sigma,\varphi)\in \mathcal{C}$, $\varphi$ provides no more information,
if we sum over all $\sigma\in \S_n$ both sides, we will get
$$\begin{array}{l}
\displaystyle \quad \sum_{\sigma\in \S_n} (-1)^{\sigma}\Delta(\sigma,\varphi,n+1)
 = \sum_{\sigma\in \S_n} (-1)^{\sigma}\big(\Delta_{\sigma(n)n}\cdots\Delta_{\sigma(1)1}\big)(\partial)\\
=\left|\begin{matrix}
x_{11} & \cdots &x_{1n}\\
\vdots & \ddots & \vdots\\
x_{n1} &\cdots & x_{nn}
\end{matrix}\right|
\cdot
\left|\begin{matrix}
\fracp{x_{11}} & \cdots &\fracp{x_{11}}\\
\vdots  & \ddots & \vdots\\
\fracp{x_{n1}} & \cdots & \fracp{x_{nn}}
\end{matrix}\right|
\end{array}$$
equals to
$$\begin{array}{l}
\displaystyle \sum_{\textrm{Capelli configuration }(\sigma,\varphi)} (-1)^{\sigma}\Delta(\sigma,\varphi,1)\\
\displaystyle = \sum_{\sigma\in \S_n} (-1)^{\sigma}
(\Delta_{\sigma(n)n}(\partial)+\delta_{\sigma(n)n}(n-1))\cdots (\Delta_{\sigma(1)1}(\partial)+\delta_{\sigma(1)1}1)\\
=\left|\begin{matrix}
D_{nn}+(n-1) & D_{n,n-1} & \cdots & D_{n1}\\
D_{n-1,n} & D_{n-1,n-1}+1 & \cdots & D_{n-1,1}\\
\vdots & \vdots & \ddots & \vdots\\
D_{1n} & D_{1,n-1} & \cdots & D_{11}
\end{matrix}\right|.
\end{array}$$
By an exchange of variables $x_1\leftrightarrow x_n, x_2\leftrightarrow x_{n-1}$, \dots, we get Cappeli's identity.

\section{The Capelli–-Cauchy–-Binet identity}

Note that, we can partially run the algorithm,
$$(-1)^{\sigma}\Delta(\sigma,\varphi,n+1)=\sum_{\Lambda^n\cdots \Lambda^{m+1}(\sigma' ,\varphi')=(\sigma,\varphi)}(-1)^{\sigma' }\Delta(\sigma' ,\varphi',m+1). $$
If we sum them up through out all $\sigma$ fixed $1,\ldots,m$,
we will get the \emph{Capelli–-Cauchy–-Binet identity}.
More exactly, for a Capelli configuration $(\sigma,\varphi)\in \mathcal{C}^{m+1}$, it is more or less clear to see from our algorithm that
$$\textrm{
$\Lambda^n\cdots \Lambda^{m+1}(\sigma)$ fix $1,\ldots,m$}\iff \begin{minipage}{0.4\linewidth}
$\sigma$ itself fix them, and
they do not appear in the image of $\varphi$.
\end{minipage}$$
If we denote the configuration with these properties by $\mathcal{C}^{m+1}_*$,
then the sum for right hand side is
$$\begin{array}{l}
\displaystyle \quad \sum_{(\sigma',\varphi')\in \mathcal{C}^m_*}(-1)^{\sigma'}\Delta(\sigma' ,\varphi',m)\\
\displaystyle =\sum_{(\sigma,\varphi)\in \mathcal{C}^m_*}(-1)^{\sigma'}
\tilde{\Delta}_{\sigma(n)n}(\partial) \cdots \tilde{\Delta}_{\sigma(m),m}(\partial) \cdot
\big(\tilde{\Delta}_{\sigma(m-1),m-1}\cdots \tilde{\Delta}_{\sigma(1)1}\big)(\partial)\\
\displaystyle = \left|\begin{matrix}
D_{nn}+(n-m) & D_{n,n-1} & \cdots & D_{n,m+1}\\
D_{n-1,n} & D_{n-1,n-1}+1 & \cdots & D_{n-1,m+1}\\
\vdots & \vdots & \ddots & \vdots\\
D_{m+1,n} & D_{1,n-1} & \cdots & D_{m+1,m+1}
\end{matrix}\right|\cdot \big(\tilde{\Delta}_{\sigma(m)m}\cdots \tilde{\Delta}_{\sigma(1)1}\big)(\partial).
\end{array}$$
The sum for left hand side is
$$\begin{array}{l}
\displaystyle \quad \sum_{\textrm{$\sigma$ fix $1,\ldots,m$}}\Delta(\sigma,\varphi,n+1)\\
\displaystyle = \sum_{\textrm{$\sigma$ fix $1,\ldots,m$}}\bigg(\tilde{\Delta}_{\sigma(n)n}\cdots \tilde{\Delta}_{\sigma(m+1),m+1}
\tilde{\Delta}_{mm}\cdots\tilde{\Delta}_{11}\bigg)(\partial)\\
\displaystyle =\sum_{\textrm{$\sigma$ fix $1,\ldots,m$}}\bigg(\tilde{\Delta}_{\sigma(n)n}\cdots \tilde{\Delta}_{\sigma(m+1),m+1}\bigg)(\partial)
\cdot\bigg(\tilde{\Delta}_{mm}\cdots \tilde{\Delta}_{11}\bigg)(\partial).
\end{array}$$
Then the classic Cauchy–-Binet argument makes sure the following theorem.

\begin{Th}[Capelli–-Cauchy–-Binet identity]If $m<n$, then
$$\begin{array}{l}
\quad \left|\begin{matrix}
D_{11}+(n-m) & D_{n,n-1} & \cdots & D_{1m}\\
D_{21} & D_{n-1,n-1}+1 & \cdots & D_{2m}\\
\vdots & \vdots & \ddots & \vdots\\
D_{mm} & D_{1,n-1} & \cdots & D_{mm}
\end{matrix}\right|\\
\displaystyle =\sum_{1\leq i_1<\cdots <i_m\leq n}
\left|\begin{matrix}
x_{1i_1} & \cdots &x_{1i_m}\\
\vdots & \ddots & \vdots\\
x_{mi_1} &\cdots & x_{mi_m}
\end{matrix}\right|
\cdot
\left|\begin{matrix}
\fracp{x_{1i_1}} & \cdots &\fracp{x_{1i_m}}\\
\vdots  & \ddots & \vdots\\
\fracp{x_{mi_1}} & \cdots & \fracp{x_{mi_m}}
\end{matrix}\right|.
\end{array}$$
\end{Th}

\section{The Turnbull's identity}
Let $x_{ij}=x_{ji}$ be a system of indeterminants. Denote
$$S_{ij}=\sum_{k=1}^n (1+\delta_{jk})x_{ik}\fracp{x_{jk}}. $$
The \emph{Turnbull's identity} claims
$$\det(S_{ij}-\delta_{ij}(n-i))=\det(x_{ij})\det\bigg((1+\delta_{jk}\fracp{x_{jk}})\bigg). $$

\bigbreak
To prove this, similarly, we define
$$\Sigma_{ij}=\sum_{k=1}^n  (1+\delta_{jk})x_{ik}y_{jk}, $$
and for $0\leq m\leq n$, and a Capelli configuration $(\sigma,\varphi)$
$$\Sigma(\sigma,\varphi,m)
=\tilde{\Sigma}_{\sigma(n)n}(\partial) \cdots \tilde{\Sigma}_{\sigma(m)m}(\partial) \cdot
\big(\tilde{\Sigma}_{\sigma(m-1),m-1}\cdots \tilde{\Sigma}_{\sigma(1)1}\big)(\partial),$$
where
$$\tilde{\Sigma}_{\sigma(i)i}= \begin{cases}
\Sigma_{\sigma(i)i}, & \textrm{if $i\neq \sigma(i)$ or $\varphi(i)=i$},\\
1, & \textrm{otherwise}.
\end{cases}$$
But now in general we do not have
$$(-1)^{\sigma}\Sigma(\sigma,\varphi,m+1)=\sum_{\Lambda^m(\sigma' ,\varphi')=(\sigma,\varphi)}(-1)^{\sigma' }\Sigma(\sigma' ,\varphi',m). $$
But when we sum up, the difference cancelled.

\begin{Lemma}Under the notations above, we have
$$\sum_{(\sigma,\varphi)\in \mathcal{C}^{m+1}}
(-1)^{\sigma}\Sigma(\sigma,\varphi,m+1)=\sum_{(\sigma,\varphi)\in \mathcal{C}^{m}}(-1)^{\sigma' }\Sigma(\sigma' ,\varphi',m). $$
\end{Lemma}

Before the proof, we firstly do some calculation.
Let $1\leq s,t,p,q\leq n$. If $t\neq  p$,  then
\def\nfracp#1{\fracp{#1}(\cdots)}%\phantom{\fracp{#1}}}
$$\begin{array}{l}
\quad \displaystyle \sum_{k,h} \bigg((1+\delta_{tk})x_{sk}\fracp{x_{tk}}\bigg)\bigg((1+\delta_{qh})x_{ph}\nfracp{x_{qh}}\bigg)\\
\displaystyle = \sum_{k,h} (1+\delta_{tk})(1+\delta_{qh})x_{sk}x_{ph}\fracp{x_{tk}}\nfracp{x_{qh}}+ \fbox{$x_{sp}$}\nfracp{x_{qt}}. %(t=h\neq p=k)
\end{array}$$
If $t=p$, then
$$\begin{array}{l}
\quad \displaystyle \sum_{k,h} \bigg((1+\delta_{tk})x_{sk}\fracp{x_{tk}}\bigg)\bigg((1+\delta_{qh})x_{ph}\nfracp{x_{qh}}\bigg)\\
\displaystyle = \sum_{k,h} (1+\delta_{tk})(1+\delta_{qh})x_{sk}x_{ph}\fracp{x_{tk}}\nfracp{x_{qh}}\\
\displaystyle \quad + \sum_{k} (1+\delta_{tk})(1+\delta_{qk})x_{sk}\nfracp{x_{qk}}\\
\displaystyle = \sum_{k,h} (1+\delta_{tk})(1+\delta_{qh})x_{sk}x_{ph}\fracp{x_{tk}}\nfracp{x_{qh}}\\
\displaystyle \quad + \sum_{k} (1+\delta_{qk})x_{sk}\nfracp{x_{qk}}+ \fbox{$x_{sp}$}\nfracp{x_{qt}}.
\end{array}$$

To simplify the notations, denote $\partial_{ij}=(1+\delta_{ij})\fracp{x_{ij}}$.

%Write $\mathcal{C}^{m+1}=C_1\sqcup C_2$, with
%$$\begin{array}{l}
%C_1=\{(\sigma,\varphi)\in \mathcal{C}^{m+1}:\sigma^{-1}(m)\geq m\}, \\
%C_2=\{(\sigma,\varphi)\in \mathcal{C}^{m+1}:\sigma^{-1}(m)< m\}.
%\end{array}$$
%So
%$$(\sigma,\varphi)\in C_i\iff \#\{(\sigma',\varphi'):\Lambda^{m}(\sigma',\varphi')=(\sigma,\varphi)\}=i. $$
%When $i< m$,
%$$\sigma'=(\sigma(i)\sigma(m))\circ \sigma\in C_1\iff \sigma\in C_1. $$
%Since if $\sigma'(i)=\sigma(m)$, then $\sigma(m)=m$.
%%Note that
%Now, $\sigma^{-1}(m)\geq m$, so
%$$\begin{array}{l}
%\quad \displaystyle \sum_{(\sigma,\varphi)\in \mathcal{C}^{m}}
%(-1)^{\sigma}\Sigma(\sigma,\varphi,m)\\
%\displaystyle = \sum_{(\sigma,\varphi)\in \mathcal{C}^{m}}\sum_{k_1,\ldots,k_n}
%(-1)^{\sigma}\big(\cdots\big)%\bigg(x_{\sigma(m+1)k_{m+1}}\partial_{x_{m+1,k_{m+1}}}\bigg)
% \bigg(x_{\sigma(m)k_m}\cdots x_{\sigma(1)k_1}\partial_{mk_m}\cdots \partial_{1k_1} \bigg)\\
%\displaystyle = \sum_{(\sigma,\varphi)\in \mathcal{C}^{m}}\sum_{k_1,\ldots,k_n}
%(-1)^{\sigma}\big(\cdots\big)%\bigg(x_{\sigma(m+1)k_{m+1}}\partial_{x_{m+1,k_{m+1}}}\bigg)
%\bigg(x_{\sigma(m)k_m}\cdots \partial_{mk_m} x_{\sigma(1)k_1}\cdots \partial_{1k_1} \bigg)\\
%\displaystyle - \sum_{(\sigma,\varphi)\in \mathcal{C}^{m}, \sigma^{-1}(m)=1}
%(-1)^{\sigma}\big(\cdots\big)
% \bigg(x_{\sigma(m)k_m}\cdots\cdots \partial_{1k_1} \bigg)\\
%\\
%\displaystyle - \sum_{(\sigma,\varphi)\in \mathcal{C}^{m}}\sum_{k_1,\ldots,k_n}
%(-1)^{\sigma}\big(\cdots\big)
%\bigg(x_{\sigma(m)\sigma(1)}\cdots  \cdots\partial_{1k_1} \bigg)\\
%\\
%\end{array}$$
$$\begin{array}{l}
\quad \displaystyle \sum_{(\sigma,\varphi)\in \mathcal{C}^{m+1}}
(-1)^{\sigma}\Sigma(\sigma,\varphi,m+1)\\
\displaystyle = \sum_{(\sigma,\varphi)\in \mathcal{C}^{m+1}}\sum_{k_1,\ldots,k_n}
(-1)^{\sigma}\big(\cdots\big)%\bigg(x_{\sigma(m+1)k_{m+1}}\partial_{x_{m+1,k_{m+1}}}\bigg)
\bigg(x_{\sigma(m)k_m}\cdots x_{\sigma(1)k_1}\partial_{mk_m}\cdots \partial_{1k_1} \bigg)\\
\displaystyle = \sum_{(\sigma,\varphi)\in \mathcal{C}^{m+1}}\sum_{k_1,\ldots,k_n}
(-1)^{\sigma}\big(\cdots\big)%\bigg(x_{\sigma(m+1)k_{m+1}}\partial_{x_{m+1,k_{m+1}}}\bigg)
\bigg(x_{\sigma(m)k_m}\cdots \partial_{mk_m}x_{\sigma(1)k_1}\cdots \partial_{1k_1} \bigg)\\
\displaystyle + \sum_{(\sigma,\varphi)\in \mathcal{C}^{m+1}, \sigma^{-1}(m)=1}\sum_{k_1,\ldots,\widehat{k_m}\ldots,k_n}
(-1)(-1)^{\sigma}\big(\cdots\big)
 \bigg(x_{\sigma(m)k_1}\cdots\cdots \partial_{1k_1} \bigg)\\
\\
\displaystyle - \sum_{(\sigma,\varphi)\in \mathcal{C}^{m}}\sum_{k_2,\ldots,\widehat{k_m},\ldots,k_n}
(-1)^{\sigma}\big(\cdots\big)
\bigg(\fbox{$x_{\sigma(m)\sigma(1)}$}\cdots  \cdots\partial_{1k_1} \bigg)
\end{array}$$
We can define an involution over $\mathcal{C}^{m+1}$. Let $(\sigma,\varphi)\in \mathcal{C}^{m+1}$, define
$(D\sigma,D\varphi)$ by resetting the value of $\sigma$ and $\varphi$ on $1$ and $m$ by
$$D\sigma(1)=\sigma(m), \quad D\sigma(m)=\sigma(1), $$
and if $D\varphi(1)=D\sigma(1)$ and $D\varphi(m)=D\sigma(m)$ if $i$ or $m$ is fixed.
This is clear a involution fixed no point, since $\mathcal{C}^{m+1}$
is exactly the figurations when the value of $\varphi$ on $\{1,\ldots,m\}\cap \Fix \sigma$ the same to $\sigma$.
But $(-1)^{\sigma}+(-1)^{D\sigma}=0$ and $x_{\sigma(1)\sigma(m)}=x_{\sigma(m)\sigma(1)}$.
So the last term of above summation cancelled.

It is not hard to see that the second term is
$$\sum_{(\sigma,\varphi)\in \mathcal{C}^{m+1}, \sigma^{-1}(m)=1}\Sigma(\hat{\sigma},\hat{\varphi},m). $$
Where $(\hat{\sigma},\hat{\varphi})$ is
the only $(\sigma,\varphi)\in \mathcal{C}^m$ other than $(\hat{\sigma},\hat{\varphi})$
with $\Lambda^m(\sigma.\varphi)=(\hat{\sigma},\hat{\varphi})$ when $\sigma^{-1}(m)<m$.
So continue this process, then we will get the expression in the lemma.

\bigbreak
Here are two remarks.
Firstly, through out the proof, we see there is no expected direct ``Cauchy–-Binet'' analogy of the Turnbull's identity.
Secondly, the Howe-–Umeda-–Kostant-–Sahi identity can be proved in the similar way.

\section{The Cayley identity}

It is an easy linear algebra exercise that $\det\left(\fracp{x_{ij}}\right)\det(x_{ij})=n!$.
%$$\left|\begin{matrix}
%\fracp{x_{11}} & \cdots &\fracp{x_{1n}}\\
%\vdots  & \ddots & \vdots\\
%\fracp{x_{n1}} & \cdots & \fracp{x_{nn}}
%\end{matrix}\right|\cdot\left|\begin{matrix}
%x_{11} & \cdots &x_{1n}\\
%\vdots & \ddots & \vdots\\
%x_{n1} &\cdots & x_{nn}
%\end{matrix}\right|=n!, \eqno{(*)}$$
%where $x_{ij}$ are now not as ``multiplying by $x_{ij}$''.
%
%Actually,
%$$D_{pq}\det(x_{ij})=\sum_{k} x_{pk}\fracp{x_{qk}}\det(x_{ij})=\det(x_{ij})|_{x_{qk}\leftarrow x_{qk}}=\delta_{pq}\det(x_{ij}). $$
%So $D_{ij}$ acts as $s$ on $\det(x_{ij})$.
%So the left hand side of Capelli's identity acts as $\det(\operatorname{diag}(s+n-1,\cdots,s)$ on $\det(x_{ij})^s$.
%Cancelling the $\det(x_{ij})$ both sides, we will get what we want.
%
It is curious what the right hand as an operator should be.
Actually, using the trick above, if we denote
$$d_{ij}=\sum_{k=1}^n \fracp{x_{ik}}x_{jk}=\sum_{k=1}^n x_{jk}\fracp{x_{ik}}+n\delta_{ij}=D_{ji}+n\delta_{ij}, $$
we will have
$$\left|\begin{matrix}
d_{11}-(n-1) & d_{12} & \cdots & d_{1n}\\
d_{21} & d_{22}-(n-2) & \cdots & d_{2n}\\
\vdots & \vdots & \ddots & \vdots\\
d_{n1} & d_{n2} & \cdots & d_{nn}
\end{matrix}\right|=
\left|\begin{matrix}
\fracp{x_{11}} & \cdots &\fracp{x_{1n}}\\
\vdots  & \ddots & \vdots\\
\fracp{x_{n1}} & \cdots & \fracp{x_{nn}}
\end{matrix}\right|\cdot\left|\begin{matrix}
x_{11} & \cdots &x_{1n}\\
\vdots & \ddots & \vdots\\
x_{n1} &\cdots & x_{nn}
\end{matrix}\right|, $$
as an operator on polynomial in $\{x_{ij}\}$.
%$$d_{ij}=\sum_{k=1}^n \fracp{x_{ik}}x_{jk}=\sum_{k=1}^n x_{jk}\fracp{x_{ik}}-\delta_{ij}=D_{ji}-\delta_{ij}, $$
Here is the sketch,
\begin{itemize}
\item Firstly change the definition of $F(\partial)$ by $H(\partial)G(x)$ as defined in the main part.
\item Then use a new $\tilde{\Delta}_{\sigma(i)i}= \begin{cases}
\Delta_{\sigma(i)i}, & \textrm{if $i\neq \sigma(i)$ or $\varphi(i)=i$},\\
-1, & \textrm{otherwise}.
\end{cases}$ to define a new $\Delta(\sigma,\varphi,m)$.
\item By a completely the same trick, we have
$$(-1)^{\sigma}\Delta(\sigma,\varphi,m+1)=\sum_{\Lambda^m(\sigma' ,\varphi')=(\sigma,\varphi)}(-1)^{\sigma' }\Delta(\sigma' ,\varphi',m). $$
\item So it follows the same argument we used before.
\end{itemize}
Note that the left hand side acts $s-1+n$ on $\det(x_{ij})^{s-1}$,
so the left hand side acts as $\det\operatorname{diag}(s,s+1,\ldots,s+n-1)$.
So we get the \emph{Cayley identity}, which claimed that
$$\left|\begin{matrix}
\fracp{x_{11}} & \cdots &\fracp{x_{1n}}\\
\vdots  & \ddots & \vdots\\
\fracp{x_{n1}} & \cdots & \fracp{x_{nn}}
\end{matrix}\right|\cdot\left|\begin{matrix}
x_{11} & \cdots &x_{1n}\\
\vdots & \ddots & \vdots\\
x_{n1} &\cdots & x_{nn}
\end{matrix}\right|^s=s(s+1)\cdots (s+n-1) \left|\begin{matrix}
x_{11} & \cdots &x_{1n}\\
\vdots & \ddots & \vdots\\
x_{n1} &\cdots & x_{nn}
\end{matrix}\right|^{s-1}.$$
Of course, this also can be get from Capelli's identity by applying $\det(x_{ij})^k$ and finally dividing itself.
%$d_{pq}$ acts as $(s-1)\delta_{pq}$ over $\det(x_{ij})^{s-1}$.

%
%\def\perm{\operatorname{perm}}
%
%\paragraph{For permutations} For a matrix $(x_{ij})$, we can denote and define
%$$\perm(x_{ij})=\sum_{\sigma\in \S_n} x_{\sigma(1)1}\cdots x_{\sigma(n)n}, $$
%similar to determinants but no signs. Then we have
%$$\perm(D_{ij}-\delta_{ij}(n-i))=\perm\big(\sum_{k=1}^n x_{ik}y_{jk}\big)\big|_{y=\partial}, $$
%by the similar argument.
%
%
%\paragraph{The transportation of Capelli's identity} If we define
%$$\overline{\Delta}(\sigma,\varphi,m)
%=\big(\tilde{\Delta}_{1\sigma(1)}\cdots \tilde{\Delta}_{m-1,\sigma(m-1)}\big)(\partial) \cdot
%\tilde{\Delta}_{m,\sigma(m)}(\partial)\cdots \tilde{\Delta}_{n,\sigma(n)}(\partial),$$
%%still with
%%$$\tilde{\Delta}_{i\sigma(i)}= \begin{cases}
%%\Delta_{i\sigma(i)}, & \textrm{if $i\neq \sigma(i)$ or $\varphi(i)=i$},\\
%%1, & \textrm{otherwise}.
%%\end{cases}$$
%Then we will get
%$$\left|\begin{matrix}
%D_{11} & D_{21} & \cdots & D_{n1}\\
%D_{12} & D_{22}+1 & \cdots & D_{n2}\\
%\vdots & \vdots & \ddots & \vdots\\
%D_{1n} & D_{2n} & \cdots & D_{nn}+(n-1)
%\end{matrix}\right|=\left|\begin{matrix}
%x_{11} & \cdots &x_{1n}\\
%\vdots & \ddots & \vdots\\
%x_{n1} &\cdots & x_{nn}
%\end{matrix}\right|
%\cdot
%\left|\begin{matrix}
%\fracp{x_{11}} & \cdots &\fracp{x_{1n}}\\
%\vdots  & \ddots & \vdots\\
%\fracp{x_{n1}} & \cdots & \fracp{x_{nn}}
%\end{matrix}\right|. $$

\vfill
\bibliographystyle{plain}
\bibliography{bibfile}

\end{document}